\documentclass[]{article}
\usepackage{amsmath, amssymb, mathrsfs, graphicx,tikz}
\usepackage{graphicx,amssymb,amsfonts,latexsym,amsmath,amsthm,times}
\usepackage{epsfig}
\usepackage{fancyhdr}
\usepackage{color}
\usepackage{cite}
\usepackage{caption}
\usepackage{authblk}
\usepackage{subcaption}
\usepackage{enumitem}
\usepackage[export]{adjustbox}
\usepackage[colorlinks,citecolor=blue,linkcolor=blue,urlcolor=blue]{hyperref}\usepackage{tikz}

\usetikzlibrary{matrix}

\newcount\colveccount
\newcommand*\colvec[1]{
        \global\colveccount#1
        \begin{pmatrix}
        \colvecnext
}
\def\colvecnext#1{
        #1
        \global\advance\colveccount-1
        \ifnum\colveccount>0
                \\
                \expandafter\colvecnext
        \else
                \end{pmatrix}
        \fi
}

\begin{document}
\title{On the Structure of Littlewood Polynomials and their Zero Sets}
\author[1]{R. Reyna}
\author[2]{S. B. Damelin}

\affil[1]{Department of Mathematics and Statistics, California State Polytechnic University Pomona. email: raphaelreyna@cpp.edu}
\affil[2]{Mathematical Reviews, The American Mathematical Society. email: damelin@umich.edu}

\maketitle

	\quad They say a picture is worth a thousand words. It is often the
    case in mathematics that a visualization of the objects we work with 
    is worth a thousand equations; visualizations allow us to observe the connections hidden behind layers of rigor and equations. In the case of Figure 1, it is worth
    67,108,856 equations, all of which are of the form $p(z) = 0$. Here, 
    all $p(z)$ are Littlewood polynomials; Figure 1 is simply the plot of polynomial roots. Specifically, Littlewood polynomials are polynomials with coefficients of either 1 or -1, and were first studied in depth by J. E. Littlewood during the 1950's. Thanks to recent advances in computing, it is now possible to compute the zeros of this class of polynomials, under some restrictions. Surprisingly, when the zeros for all Littlewood polynomials, up to a certain degree, are plotted, (see Figure 1 and Figure 2), the zero sets are locally fractal-like and exhibit high degrees of symmetry. \newline \par Motivated by the complex symmetry of the zero set of Littlewood polynomials, we examine the symmetry of the set of Littlewood polynomials themselves.
    In Figure 2 we see that the distribution of the zero sets are locally fractal-like. Thus, the complexity of the set grows exponentially on a local scale. The computational complexity of such sets is to be expected, since computing this set means solving the root finding problem for potentially millions of polynomials. We use the relationship between the symmetry in the zero set of Littlewood polynomials and in the symmetry in their coefficients to reduce the computational complexity of the zero set and to give an alternate characterization of the set. Namely, we will study connections between the group structure of the coefficients of Littlewood polynomials, the distributions of the zeros and a fractal set called a generalized dragon set which will be defined along the way.
    
    \begin{figure}[h!]
 \centering
 \includegraphics[scale=3.5]{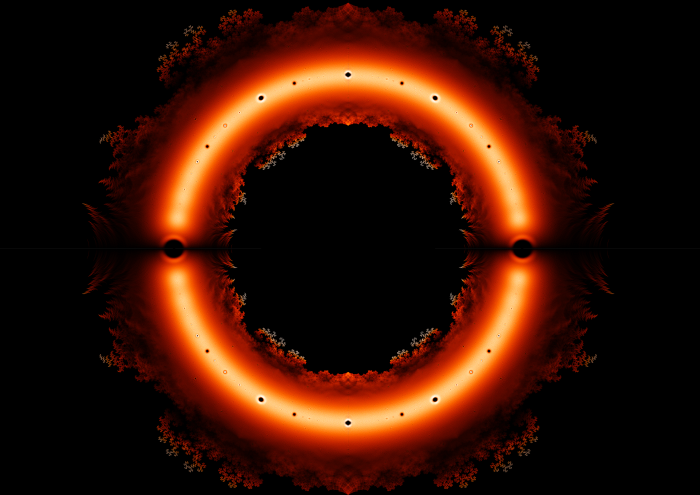}
 \caption{Roots of Littlewood polynomials up to degree 24 in $[-2,2]\times[-1.5,1.5] \subset \mathbb{C}$. Coloration denotes density.\newline Image credit: S. Derbyshire}
\end{figure}

    \begin{figure}[h!]
 \centering
 \includegraphics[scale=0.1]{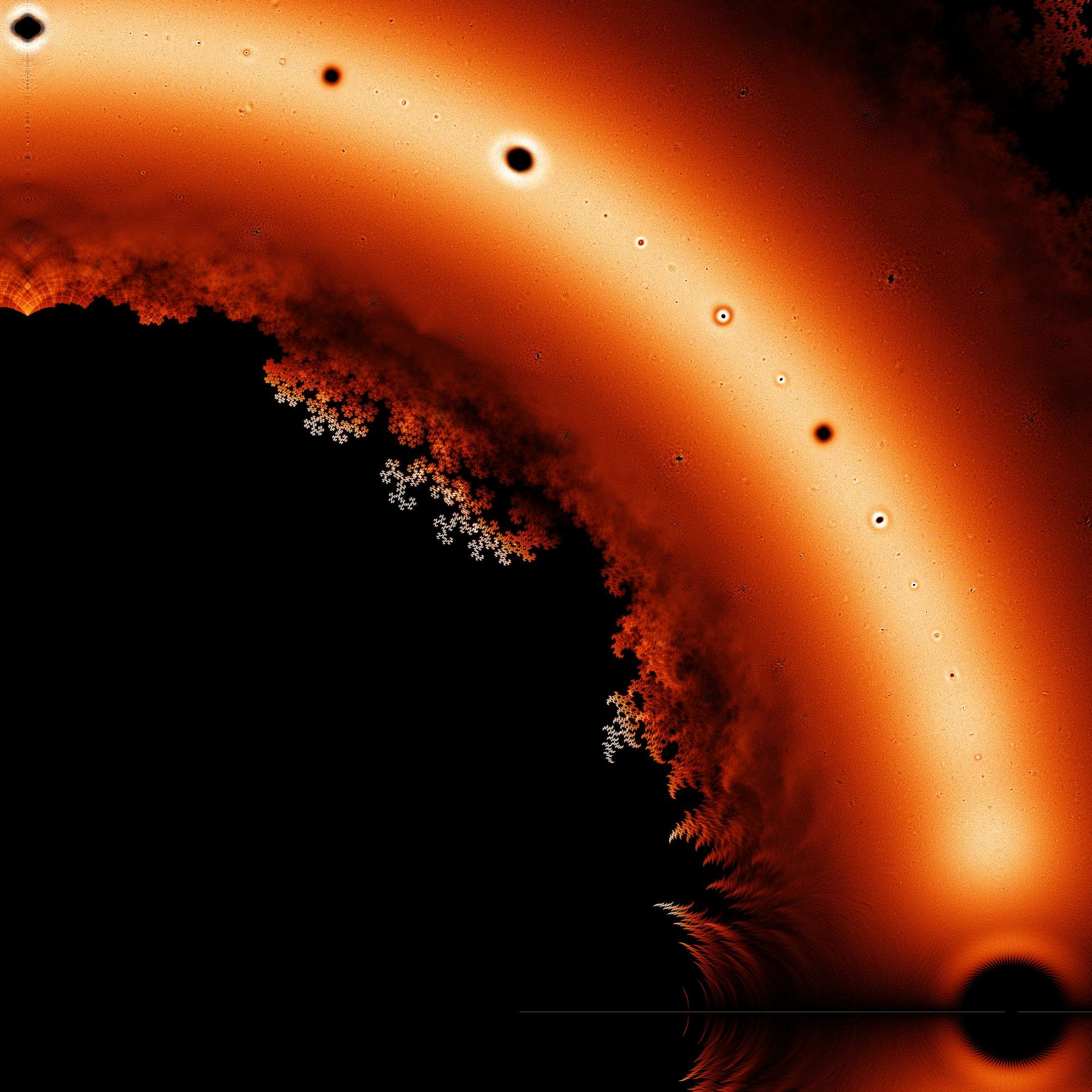}
 \caption{Close up of Figure 1, first quadrant of $\mathbb{C}$. \newline Image credit: S. Derbyshire}
\end{figure}

\par Before we jump into exploring the symmetries of fractals and fractal-like objects, we should be clear on what a fractal is. While a universally accepted definition of a fractal set has yet to be established in the literature, in the famous words of US Supreme Court Justice, Potter Stewart, "I know it when I see it". This very unmathematical and fuzzy manner of identifying fractals and fractal-like objects is perhaps due to fractal geometry being a very new branch of mathematics. Fortunately, if one is a bit more selective with the type of fractals one wishes to study, it is possible to provide a solid definition. For our purposes, we shall concern ourselves with \textit{iterated function systems}, usually written \textit{IFS}. An IFS is a finite set of contraction mappings on a complete metric space (for this paper the reader may just assume that by \textit{complete metric space} we mean $\mathbb{R}^2$ or $\mathbb{C}$). The way one constructs a fractal from an IFS is by iterated applications of contraction mappings to a compact set. If we let 
$\{f_j:X\rightarrow X\ |\ j=0,1,2,\hdots,d\}$, be a set of contraction mappings on a compact space endowed with the Hausdorff metric, $X$, (here, the reader may assume $X$ is roughly $\mathbb{R}^2$ or $\mathbb{C}$ along with a notion of distance between subsets) then, by the Banach fixed-point theorem \cite{KF}, there exists an attractor set 
$$S=\lim_{n\rightarrow \infty}S_n,$$ where 
$$S_{n+1} = \bigcup_{j=1}^d f_j(S_n),$$
which will be called a fractal.\newline
\quad In order to solidify the notion of an IFS, we will show the construction of a well known fractal set, the Sierpinski Triangle. We let $X =\mathbb{C}$, $f_1(x) = \frac{1}{2}x$, $f_2(x) = \frac{1}{2}x +\frac{1}{2}$, and $f_3(x) = \frac{1}{2}x+(\frac{1}{4}+\frac{\sqrt{3}}{4}i)$ be our contraction mappings, and $S_0$  be an equilateral triangle region in $\mathbb{C}$ with vertices at $0, 1$, and $\frac{1}{2}+\frac{\sqrt{3}}{2}i$. Then, as outlined above, the set $S_1 = f_1(S_0) \cup f_2(S_0) \cup f_2(S_0)$. Thus, $S_1$ is the union of three maps of $S_0$. The plots below visually show this process for up to $n =3$.\newline\ \newline
\includegraphics[scale = 0.27]{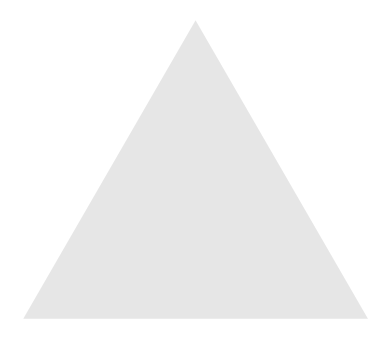}, \includegraphics[scale = 0.27]{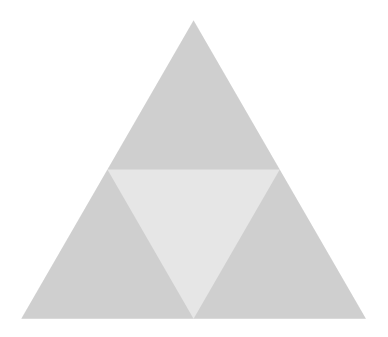}, \includegraphics[scale = 0.27]{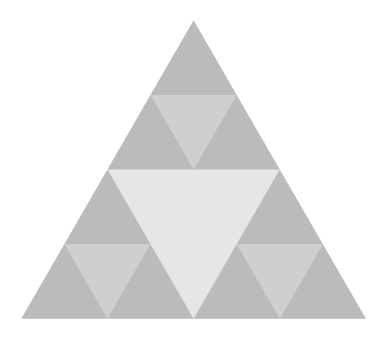},\includegraphics[scale = 0.27]{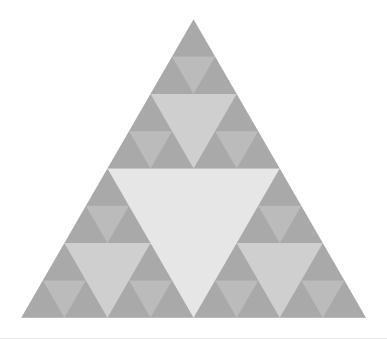}.\newline
\par In the first plot we see the region $S_0$. In the second plot, we see $S_0 \cup S_1$ with $S_1$ being a darker shade. Similarly, the darker regions of the third and fourth plot are $S_2$ and $S_3$ respectively. In short, what this process is doing is taking $S_{n-1}$, scales it down, makes three copies, and then applies different translations to each copy to produce $S_n$.\newline

\quad Another excellent, and much more relevant example of an IFS is the Littlewood polynomials themselves! Recall that the Littlewood polynomials are polynomials with coefficients of $\pm 1$. With that in mind, let us see how these polynomials may be constructed through function iteration. Let $f_+(x) = 1+zx$ and $f_-(x) = 1-zx$ be contraction mappings from $\mathbb{C}$ to $\mathbb{C}$. If we start to iterate these two contractions we get something like
$$f_+(f_-(f_+(f_+(0)))) = 1-z-z^2+z^3.$$ 
Notice that no matter how we compose $f_{-}$ and $f_{+}$, we shall always get a polynomial with coefficients of $\pm 1$. However, the set of Littlewood polynomials of some fixed degree $d$, which we shall denote as 
$$\mathcal{L}_d = \left\{p(z) = \sum_{j=0}^{d}a_jz^j\ |\ a_j \in \{-1,1\}\right\},$$
includes polynomials with $-1$ as a constant term. Therefore, by iterating $f_{+}$ and $f_{-}$ we only get half of the Littlewood polynomials, namely, the ones with a constant term of $1$. Fortunately, since we are concerned with their roots, and the IFS construction of $\mathcal{L}_d$ is only off by a factor of $-1$, this IFS construction will suffice. 
In our first example of an IFS, our final object that we ended up with was some geometric set in $\mathbb{R}^2$, but in this second example, we ended up with a set of functions. This is perfectly okay! It is exactly this fractal of functions that can then give us a "geometric" fractal. The type of fractal object that one gets from this collection of functions is what is known as a \textit{generalized dragon set}. Therefore, a generalized dragon set will be the set of all polynomials in $\mathcal{L}_d$ with a constant term of 1, evaluated at some $z \in \mathbb{C}$. In Figure 3a and Figure 3b below, we illustrate this fact below with a plot of 
$$\mathcal{L}_{13}(z_0) = \{p(z_0)\ |\ p \in \mathcal{L}_{13}\},$$ for fixed $z_0 \in \mathbb{C}$ and of its corresponding Generalized Dragon Set.

\begin{center}
	\includegraphics[scale=0.75]{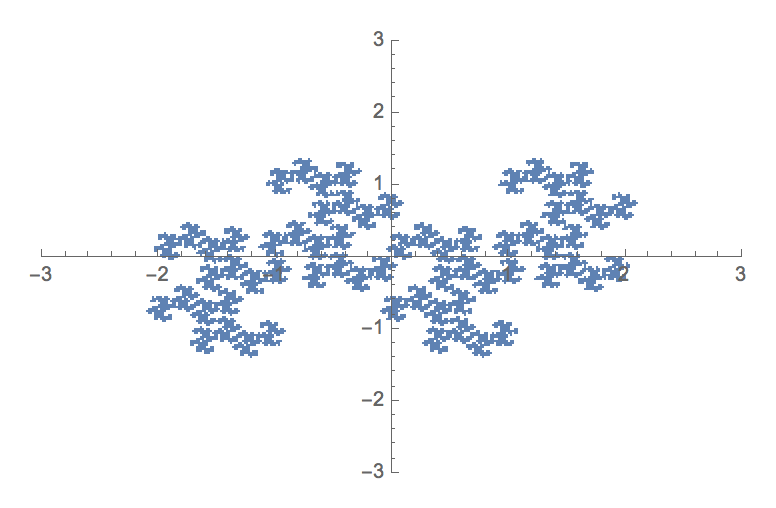}
\end{center}
\begin{center}
	\textit{\textbf{Figure 3a} All Littlewood polynomials of degree 13 evaluated at $z_0=0.48+0.45i$}
\end{center}

\begin{center}
	\includegraphics[scale=0.75]{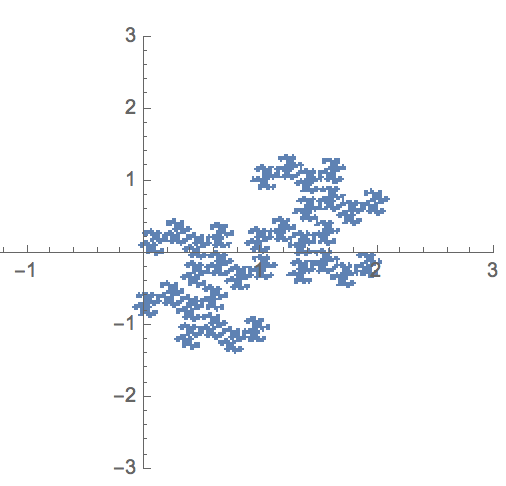}
\end{center}
\begin{center}
	\textit{\textbf{Figure 3b} 14 iterations of generalized dragon set with parameter $z_0=0.48+0.45i$}
\end{center}

Now that we have seen the connection between Littlewood polynomials and generalized dragon sets, we will continue our study of the roots of Littlewood polynomials via iterated function systems. Not only does this connection outline where the locally fractal-like structure of the roots comes from but it also allows us to exploit the high degree of symmetry in the set of polynomials to approximate the roots,
which we shall denote as
 $$D_d = \left\{z \in \mathbb{C}\ |\ p(z) = 0,\ \forall p \in \mathcal{L}_d\right\}.$$ A visualization of such a set can be found in Figure 1, as it is $D_{24}$. Naturally, we shall use groups, to capture this symmetry. Thus, we wish to find an appropriate group structure on $\mathcal{L}_d$\newline 

 Two obvious choices for the binary operation on $\mathcal{L}_d$ would be either the ordinary polynomial product, or polynomial addition. However, we immediatly see that neither of those would preserve the requirement that the coefficients be $\pm 1$. Instead, we use a matrix product from Linear Algebra, the Hadamard product, to induce a group structure to the set $\mathcal{L}$. \newline
The {\it Hadamard product}, which we shall denote as $\circ$, is defined to be component-wise multiplication between two matrices of the same size.\newline For example, the Hadamard product of two row vectors is given as $(a_1,a_2)\circ(b_1,b_2) = (a_1*b_1, a_2*b_2)$, where $*$ denotes regular multiplication of scalars.\newline
Equipping $\mathcal{L}_d$ with $\circ$ gives an Abelian group (that is, a group where the operation is commutative) with $2^{d+1}$ elements. The closure, associativity, and commutativity of $\circ$ all come from the closure, associativity, and commutativity of $*$ over the set $\{-1,1\}$. Furthermore, the identity element will be the polynomial of degree $d$ with all coefficients being 1, and every Littlewood polynomial is its own inverse in this group. The reader may have noticed by now that the group $(\mathcal{L}_d,\circ)$ is naturally isomorphic to $\bigoplus_{j=0}^{d}\mathbb{Z}_2$, where $\bigoplus$ denotes the direct sum of groups, and $\mathbb{Z}_2$ is the group of integers mod 2 under multiplication mod 2. Essentially, we pick up a copy of $\mathbb{Z}_2$ for every coefficient in a Littlewood polynomial.

It is known that $\bigoplus_{j=0}^{d}\mathbb{Z}_2$ is finitely generated, and so $\mathcal{L}_d$ must also be finitely generated by some subset $S_d \subset \mathcal{L}_d$, with $d+1$ elements. For our purpose, we shall only consider $S_d$ to be the set of Littlewood polynomials with only one negative coefficient. That is, $$S_d =\left\{ \colvec{4}{-1}{\ 1}{\ \vdots}{\ 1}, \colvec{4}{\ 1}{-1}{\ \vdots}{\ 1}, \ \cdots\ , 
	\colvec{4}{\ 1}{\ 1}{\ \vdots}{-1} \right\}.$$\newline

Now that a group structure on $\mathcal{L}_d$ has been established, the existence of a group 
generating set will allow for a coordinate system to be established on $\mathcal{L}_d(z) = \{p(z)\ : p \in \mathcal{L}_d\}$. Ultimately, the generating set for $\mathcal{L}_d$ will allow us to factor any Littlewood polynomial over the Hadamard product into a product of polynomials in the generating set. Note that this will not be factoring in the usual sense, i.e. over the canonical polynomial product (that is, multiplication of polynomials in the usual sense using the distributive property). To get a sense of polynomial multiplication in the sense of Hadamard, let us compare it with the usual product of polynomials. For example, let $p(x) = -x^2+x+1$, and $q(x) = -x^2-x+1$. Then, the Hadamard product of $p(x)$ and $q(x)$ will be: 
$$p(x)\circ q(x) = x^2-x+1.$$ Compare this with the usual way of multiplying polynomials: 
$$p(x)q(x) = x^4 - 3x^2 +1.$$
It will turn out to be the case that by factoring any polynomial in $\mathcal{L}_d$ in the sense of Hadamard, we shall be able to uniquely write any point in $\mathcal{L}_d(z_0)$ as a binary linear combination of a much smaller subset of $\mathcal{L}_d(z_0)$. That is, we shall be able to compute the evaluation of a large number of polynomials without actually evaluating too many of them. \newline The following image highlights the points in $\mathcal{L}_{13}(0.48+0.45i)$ which are the images of a polynomial in $S_{13}$ in red and the point which is the image of $e_{13}$ in green. In short, we will show that every point in blue can be expressed as a binary sum of points in red, plus a correcting term which will be a scalar multiple of the point in green. Furthermore, we will show how if a scalar multiple of the point in green can be expressed as a linear combination over $\{0,1\}$ of the points in red, then the parameter $z$ used in the construction of this image is a root of some Littlewood polynomial.

\begin{center}
	\includegraphics[scale=0.75]{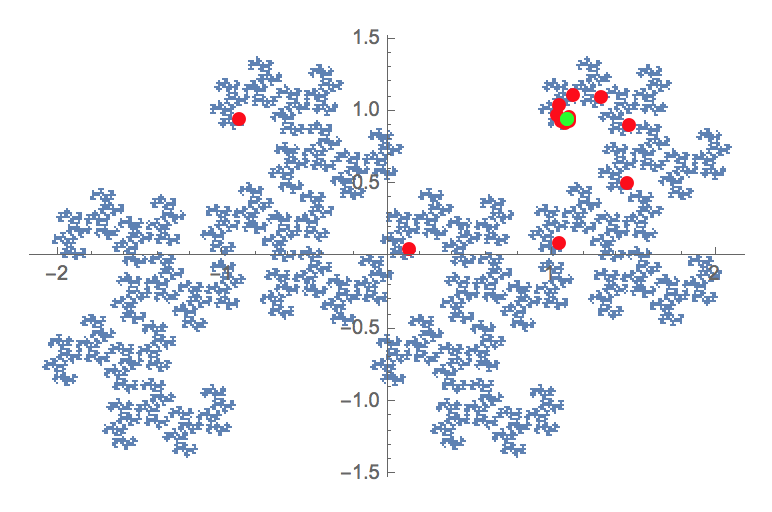}
\end{center}
\begin{center}
	\textbf{Figure 4} \ $\mathcal{L}_{13}(0.48+0.45i)$ \textit{with} $S_{13}(0.48+0.45i)$ \textit{in green, and} 
	$e_{13}(0.48+0.45i)$ \textit{in red.}
\end{center}

\subsection*{Factoring Function:}

Since the group $\mathcal{L}_d$ is generated by some subset, it is only natural to ask how do Littlewood polynomials factor over this generating set, and how does polynomial evaluation interact with the process of factoring in the sense of Hadamard. In order to gain some insight, we need to define a factoring function. Let $p \in \mathcal{L}_d$. Then, define a factor function, $\sigma(p)$, to be a subset of the generating set, $S_d$, of $\mathcal{L}_d$ such that the product of its elements is $p$. In short, since $\mathcal{L}_d$ is generated by some subset, the factoring function $\sigma$ tells us how to write any polynomial in $\mathcal{L}_d$ as the Hadamard product of Littlewood polynomials in the generating set. Furthermore, we define a measure function $\nu:\mathcal{L}_d \rightarrow \mathbb{N}$ on $\mathcal{L}_d$ such that $\nu(p)$ is the number of negative coefficients of $p$. That is,  $\nu(p) = |\sigma(p)|$, where $|\sigma(p)|$ is the number of elements in the set $\sigma(p)$.\newline
For an example of $\sigma$ and $\nu$, consider the Littlewood polynomial $x^3-x^2-x+1 \in \mathcal{L}_3$. Then, $\sigma(x^3-x^2-x+1) = \{x^3-x^2+x+1,x^3+x^2-x+1\}$, and $\nu(x^3-x^2-x+1) = 2$.

The following lemma and theorem give a way of taking the generating set of $\mathcal{L}_d$ and using it to write any point in $\mathcal{L}_d(z_0)$ as a linear combination of points in $S_d(z_0)$. Essentially, they will summarize how polynomial evaluation and the Hadamard product of polynomials interact with one another.\bigskip

{\bf Lemma 1}\,
	Let $S_d$ be the generating set of $\mathcal{L}_d$ with $d\geq 2$ and $z \in \mathbb{C}$. Then,
	$$e(z) = \frac{1}{d-1}\sum_{g \in S_d}g(z).$$
\medskip 

\begin{proof}
	We begin by observing that $e(z) = \sum_{k=0}^{d}z^k$. Furthermore,
	$$\sum_{g \in S_d}g(z) = \sum_{j=1}^{d+1}\sum_{i=0}^{d}(1-2\delta_{j}^{i})z^i=
	(d-1)\sum_{i=0}^{d}z^i.$$
	Thus, 
	$$e(z) = \frac{1}{d-1}\sum_{g \in S_d}g(z).$$
\end{proof}\bigskip

As an example of Lemma 1, we let $d = 3$ and $z_0 = 0.5+0.5i$. Then by Lemma 1, 
$$e(z_0) = \frac{1}{2}((-z_0^3+z_0^2+z_0+1)+(z_0^3-z_0^2+z_0+1) + (z_0^3+z_0^2-z_0+1) + (z_0^3+z_0^2+z_0-1))$$
$$e(0.5+0.5i) = \frac{1}{2}(1.75+0.75i + 1.25+0.25i + 0.25+0.25i-0.75+1.25i)$$
$$e(0.5+0.5i) = 1.25+1.25i.$$
In short, Lemma 1 gives a way of evaluating the identity element in $\mathcal{L}_d$ in an indirect way.\newline
The next theorem then allows us to do the same, evaluate arbitrary polynomials in the group $\mathcal{L}_d$ in an indirect manner.

{\bf Theorem 2}\, Let $p\in \mathcal{L}$. If $\circ$ is the Hadamard product then, for all $z \in \mathbb{C}$,
	$$p(z) =(1-\nu(p)) e(z) + \sum_{g \in \sigma(p)} g(z).$$
\medskip

\begin{proof}
	Let $p \in \mathcal{L}$ such that $p = g_1 \circ g_2 \circ \cdots \circ g_k$ where 
	$g_i \in S_d$ and $k = \nu(p)$. Also, let $p^j$ and $g^j$ denote the $j^{th}$ component of $p$ and $g$ respectively. Then, if $p^j = 1$ it follows 
	that $g^j_i = 1$ for all $1\leq i \leq k$. Thus,
	$$\left(\sum_{i=1}^k g^j_i\right) - (k-1) = k - (k-1) = 1 = p^j.$$
	Now if $p^j = -1$, there exists only one $g_i$ such that $g^j_i = -1$. Thus, 
	$$\left(\sum_{i=1}^kg^j_i\right) - (k-1) = (k-2)-(k-1) = -1 = p^j.$$
	Therefore, in both cases, 
	$$p^j = \left(\sum_{i=1}^{k}g^j_i\right)-(k-1).$$
	Thus, we conclude that 
	$$p = \left(\sum_{i=1}^kg_i\right)-(k-1)e.$$
    Thus, since the polynomials $p,g$ and $e$, are all in $\mathbb{C}[x]$, it follows that for the polynomial functions, one has
    $$p(z) =(1-\nu(p)) e(z) + \sum_{g \in \sigma(p)} g(z).$$
\end{proof}

As an example of Theorem 2, we let $p \in \mathcal{L}_3$ such that $p(x) = x^3-x^2-x+1$, and $z_0 = 0.5+0.5i$. Then 
$$p(z_0) = (1-2)e(z_0)+(z_0^3-z_0^2+z_0+1)+(z_0^3+z_0^2-z_0+1)$$
$$p(0.5+0.5i) = -1.25-1.25i + 1.25+0.25i + 0.25+0.25i$$
$$p(0.5+0.5i) = 0.25-0.75i.$$
It has now been shown how the group structure of $\mathcal{L}_d$, the group of Littlewood polynomials of degree $d$ under the Hadamard product, can be used to establish an addressing system on Generalized Dragon Sets, which we have already shown to be what $D_d$ looks like on a local scale. That is, any point in the Generalized Dragon set, $\mathcal{L}_d(z)$, can be expressed as a linear combination of the elements of the generating set, $S(z)$ and the identity, $e(z)$, of $\mathcal{L}_d$. This addressing system, taken together with the fact that $D_d$ locally looks like $\mathcal{L}_d(z_0)$ for some $z_0 \in \mathbb{C}$, gives us a sort of local coordinate system on $D_d$.\newline
This then allows us to essentially slice up $\mathcal{L}_d(z_0)$ into smaller, more manageable partitions. In mathematics, mathematical objects with an underlying set, such as our groups of Littlewood polynomials, are partitioned by using equivalence classes. Formally, we begin by observing that $1\leq\nu(p)\leq d$, with $\nu(p) \in \mathbb{N}$. Let $N_d(n) = \{ p \in \mathcal{L}_d\ :\ \nu(p) = n\}$ be an equivalence class given by the partition induced by $\nu$. Then, $|N_d(n)| = {d \choose \nu(p)}$, and $\bigcup_{n=0}^d N_d(n) = \mathcal{L}_d$. \newline To provide a visualization of such a partition on a dragon set, consider the plots below, which are the plots of the evaluations of polynomials in $N_{14}(n)$ at $z =0.48+0.45i $, for $1 \leq n \leq d$. Evaluations of polynomials in $S_{14}(z)$ are shown in red.\newpage

\begin{figure}[h]
 
\begin{subfigure}{0.5\textwidth}
\includegraphics[width=0.9\linewidth, height=5cm]{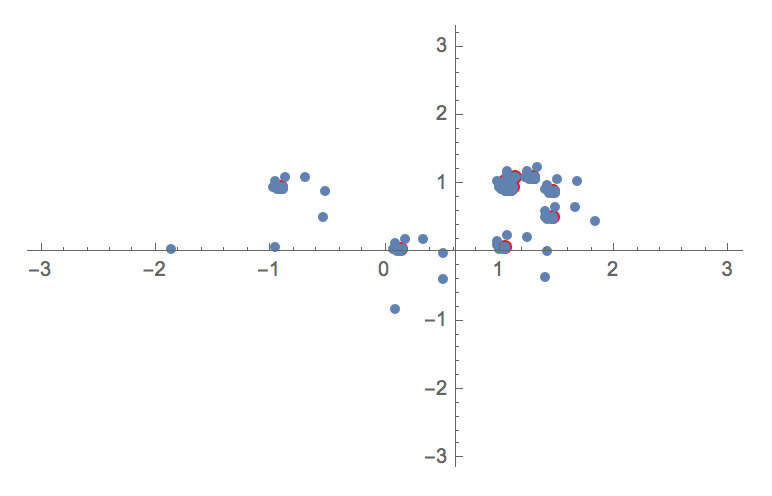} 
\caption{\textit{Evaluation at 0.48+0.45i of Littlewood polynomials of degree 14 with }$\nu(p) = 2$.}
\label{fig:subim11}
\end{subfigure}
\begin{subfigure}{0.5\textwidth}
\includegraphics[width=0.9\linewidth, height=5cm]{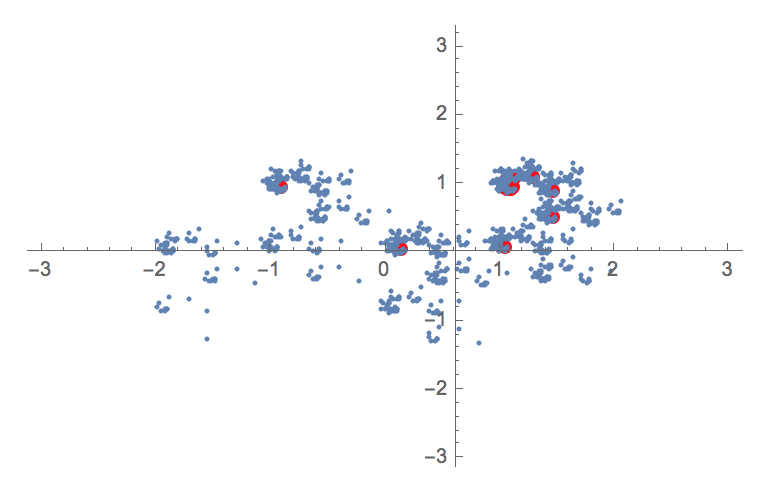}
\caption{\textit{Evaluation at 0.48+0.45i of Littlewood polynomials of degree 14 with }$\nu(p) = 4$.}
\label{fig:subim12}
\end{subfigure}
 
\caption{}
\label{fig:image1}
\end{figure}

\begin{figure}[h!]
 
\begin{subfigure}{0.5\textwidth}
\includegraphics[width=0.9\linewidth, height=5cm]{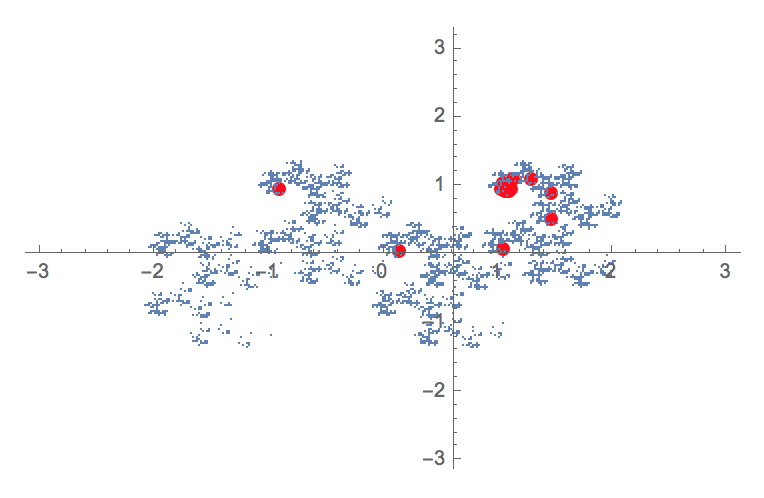} 
\caption{\textit{Evaluation at 0.48+0.45i of Littlewood polynomials of degree 14 with }$\nu(p) = 6$.}
\label{fig:subim21}
\end{subfigure}
\begin{subfigure}{0.5\textwidth}
\includegraphics[width=0.9\linewidth, height=5cm]{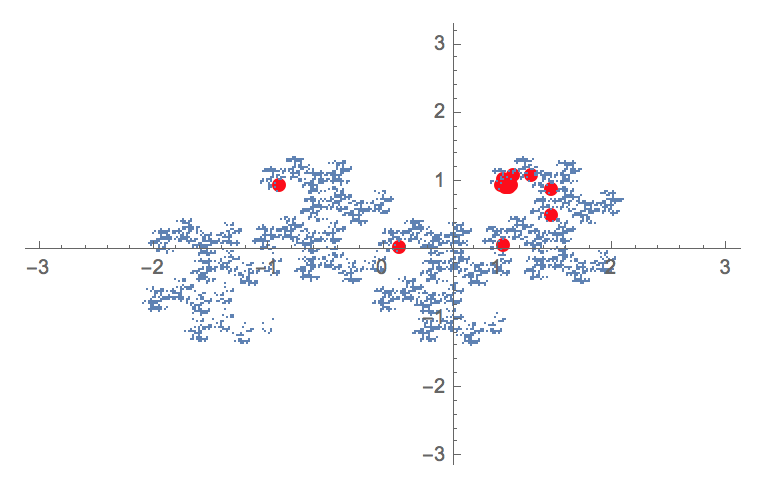}
\caption{\textit{Evaluation at 0.48+0.45i of Littlewood polynomials of degree 14 with }$\nu(p) = 7$.}
\label{fig:subim22}
\end{subfigure}
 
\caption{}
\label{fig:image2}
\end{figure}

\begin{figure}[h!]
 
\begin{subfigure}{0.5\textwidth}
\includegraphics[width=0.9\linewidth, height=5cm]{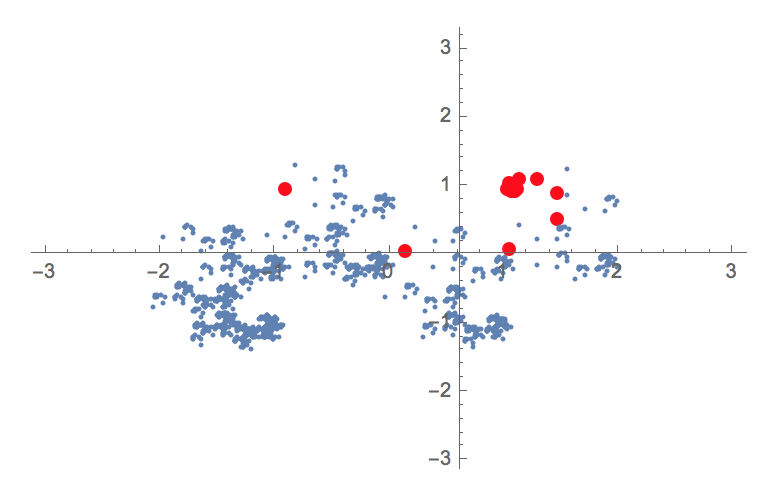} 
\caption{\textit{Evaluation at 0.48+0.45i of Littlewood polynomials of degree 14 with }$\nu(p) = 11$.}
\label{fig:subim31}
\end{subfigure}
\begin{subfigure}{0.5\textwidth}
\includegraphics[width=0.9\linewidth, height=5cm]{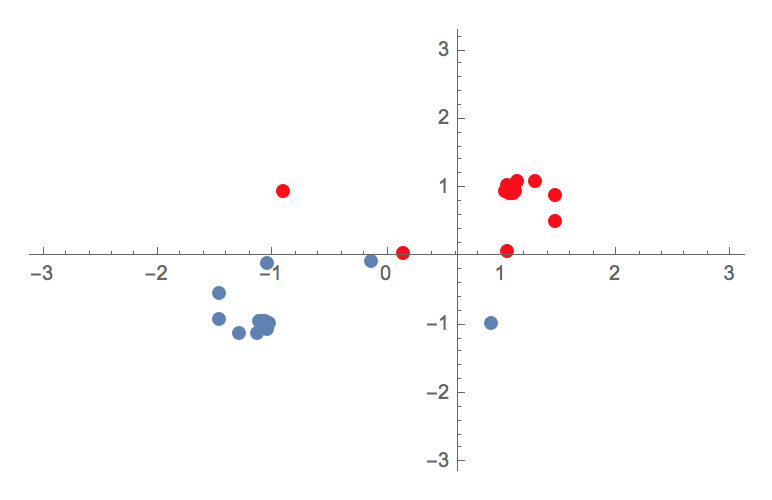}
\caption{\textit{Evaluation at 0.48+0.45i of Littlewood polynomials of degree 14 with }$\nu(p) = 14$.}
\label{fig:subim32}
\end{subfigure}
 
\caption{}
\label{fig:image3}
\end{figure}\newpage

From these plots, we observe that $\mathcal{L}_d(z)$ can be approximated by evaluating polynomials in $N_d(k) \subset \mathcal{L}_d$, where $k$ is some fixed constant in $\mathbb{N}$. When $k$ is chosen to be small, the number of polynomial evaluations is large and when $k$ is chosen to be closer to $d$ the number of polynomial evaluation is small. We observe in the above plots that the evaluations of polynomials in $N_d(k)$ is a reflection of the evaluation of polynomials with $N_d(d-k)$. Thus, an approximation of this type will require at most ${d \choose \lceil\frac{d}{2}\rceil}$ polynomial evaluations. Indeed, we observe in the above plots, for $d=14$,  choosing $k=7$ provides the best approximation. Therefore, we observe that even the ``largest" approximation of this type will be a set with ${d \choose \lceil\frac{d}{2}\rceil}$ elements as opposed to $\mathcal{L}_d(z)$ which will have $2^{d+1}$ elements. For a comparison, the following plot shows the difference in cardinalities between $\mathcal{L}_d(z_0)$ and the best approximation as shown in the figures above.

\begin{figure}[h!]
 \centering
 \includegraphics[scale=0.3]{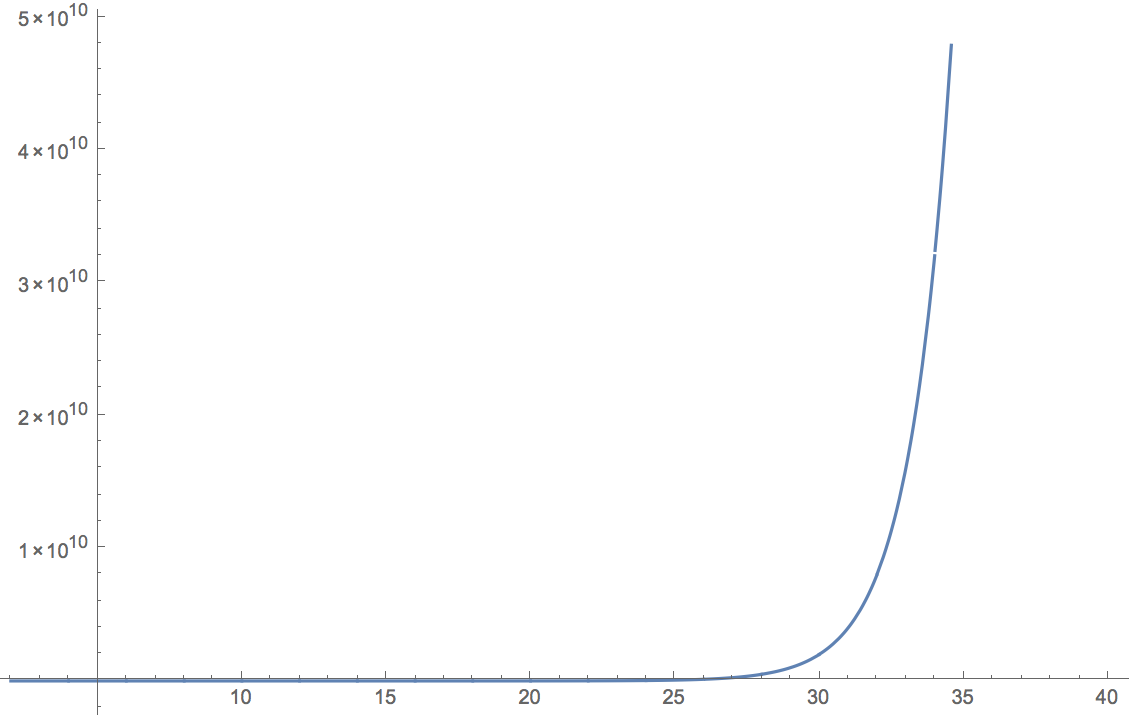}
 \caption{}
\end{figure}
    
 We see that even the ``largest" approximation of $\mathcal{L}_d(z)$ is much smaller and cheaper to compute than $\mathcal{L}_d(z)$ and becomes more efficient as the degree of the polynomials increases. Since our main interest is in determining when some given $z \in \mathbb{C}$ is a root of some polynomial in $\mathcal{L}_d$, we observe that it would be significantly cheaper to check if the ``largest" approximation intersects some $\varepsilon$-ball about 0.\newline

In conclusion, by looking at the group structure of the set of Littlewood polynomials for some fixed degree $d$, and under the appropriate group operation, we may build easily computable approximations of its zero set. Computing the zero set of $\mathcal{L}_d$ involves finding roots of a very large number of high degree polynomials, which, while it is a straighforward procedure, it is computationally very expensive. Luckily, the zero set of $\mathcal{L}_d$ has enough symmetry to allow for a much cheaper approximation.

\end{document}